\documentclass[a4paper,10pt]{article}
\usepackage{url}
\usepackage{hyperref}
\usepackage{verbatim}
\title{On the degree of the polynomial defining a planar algebraic curves of constant width}
\author{Magali Bardet\thanks{Universit\'e de Rouen, Laboratoire LITIS, \'Equipe C\&A, avenue de l'Universit\'e, F-76801 Saint-\'Etienne-du-Rouvray Cedex, (magali.bardet@univ-rouen.fr).} \and T\'erence Bayen\thanks{Universit\'e Montpellier 2, Laboratoire I3M, CC 051, Place Eug\`ene Bataillon, 34095 Montpellier cedex 5, France (terence.bayen@univ-montp2.fr).}}
\date{\today}
\usepackage[francais,english]{babel}
\usepackage[latin1]{inputenc}
\usepackage{dsfont}
\usepackage{graphicx}
\usepackage{color}
\usepackage{epsfig}
\usepackage{amsthm,amssymb,amsbsy,amsmath,amsfonts,amssymb,amscd,mathrsfs}

\newtheorem{theo}{Theorem}[section]

\newtheorem{lem}{Lemma}[section]
\newtheorem{conj}{Conjecture}[section]
\newtheorem{prop}{Proposition}[section]
\newtheorem{defi}{Definition}[section]
\newtheorem{coro}{Corollary}[section]
\newcommand{\R}{\mathbb{R}}
\newcommand{\ov}{\overline}

\newcommand{\vphi}{\varphi}
\numberwithin{equation}{section}

\newcommand\be{\begin{equation}}
\newcommand\ee{\end{equation}}
\begin{document}
\maketitle
\begin{abstract}
In this paper, we consider a family of closed planar algebraic curves $\mathcal{C}$
which are given in parametrization form via a trigonometric polynomial $p$.
When $\mathcal{C}$ is the boundary of a compact convex set, the polynomial $p$ represents the support function of this set.
Our aim is to examine properties of the degree of the defining polynomial of this family of curves in terms of the degree of $p$.
Thanks to the theory of elimination, we compute the total degree and
the partial degrees of this polynomial, and we solve in addition a question raised by Rabinowitz in \cite{Rabi}
on the lowest degree polynomial whose graph is a non-circular curve of constant width.
Computations of partial degrees of the defining polynomial of algebraic surfaces of constant width are also provided in the same way.
\end{abstract}
\section{Introduction}
We consider the set $\mathcal{E}$ of planar algebraic curves $\mathcal{C}$ defined by:
\begin{equation}{\label{def-curve-C}}
\theta \in [-\pi,\pi] \longmapsto \gamma(\theta):=
\begin{cases}
x(\theta) = p(\theta) \cos \theta - p'(\theta) \sin \theta,\\
y(\theta) = p(\theta) \sin \theta + p'(\theta) \cos \theta,
\end{cases}
\end{equation}
where $p$ is a trigonometric polynomial of degree $N \geq 1$. The previous
parametrization is standard in convex analysis to describe the boundary of a
non-empty compact convex set $K$ of $\mathbb{R}^2$ via its support function $p$, see \cite{Bon,schneider}.
More generally, given a convex set $K \subset \R^n$, $n \geq 2$, one can give a similar parametrization of its boundary via its support function $h_K$.
When $p$ is a periodic function (of period $2\pi$), it is well-known that (\ref{def-curve-C}) defines the boundary of a planar convex set $K$
if and only if $p+p'' \geq 0$ in the distribution sense, see \cite{Groemer,vitale}, and 
$\mathcal{C}$ coincides with the inverse Gauss-mapping of $K$, see \cite{Bayen4}.
It follows that the family $\mathcal{E}$ contains the class of non-empty compact convex sets of $\mathbb{R}^2$ whose support function is
a trigonometric polynomial.

The study of algebraic curves defined by (\ref{def-curve-C}) has several applications in
optimization, see e.g. reference \cite{did-ter} for a reformulation of convexity constraint by semi-definite programming
for several shape optimization problems. In \cite{zbynek}, polynomial support functions are also used for geometric design.
An interesting subclass of $\mathcal{E}$ is the set of all planar
algebraic curves of constant width. A curve of constant width $\alpha$ defines the boundary of a convex set and is given by (\ref{def-curve-C})
under the additional constraint $p(\theta)+p(\theta+\pi)=\alpha$,
for all $\theta \in [0,2\pi]$. These curves have been widely studied
since the 19th century and have applications in mechanics.
Geometrical properties of these curves can be found in the literature,
see e.g. \cite{BayenJBHU,Groemer,Leb1,Meiss1,Meiss2,Reul,Yag}.

In this work, we study properties of algebraic curves of constant,
and in particular the question of implicitization of these curves.
Rabinowitz, see \cite{Rabi}, has raised the question of finding a non-circular planar constant width curve
whose defining polynomial is of lowest degree. By considering constant width curves $\mathcal{C}$ with trigonometric
polynomial support function of degree $3$, he obtains an implicit representation of these curves by a polynomial of total degree $8$,
and he conjectures that this sub-class is of lowest degree (besides the circle of degree 2).

In the present work, we compute exactly the total degree of the
defining polynomial of any curve of $\mathcal{E}$. As a consequence, we obtain that it is impossible to obtain a 
polynomial of degree less than 8 defining a non-circular constant width curve, which answers to the question raised by Rabinowitz.
Our main result is the following (see Theorem \ref{main}): if $p$ is of degree $N \geq 2$ and has only odd coefficients, then
the total degree of the defining polynomial $f\in\mathbb{R}[x,y]$ of $\mathcal{C}$ is $N+1$, whereas if $p$ has
an even non-zero coefficient, then $f$ is of total degree $2N+2$.

The proof of this result relies on properties of the implicitization of planar curves which can be found for instance in \cite{sendra1}.
The defining polynomial of a curve of $\mathcal{E}$ is obtained as follows.
We first compute the number of times $\gamma$ traces a curve $\mathcal{C}$ in $\mathcal{E}$,
and we consider a rational parametrization of
$\mathcal{C}$ of the form $\left(\frac{\chi_{1,1}(t)}{\chi_{1,2}(t)},\frac{\chi_{2,1}(t)}{\chi_{2,2}(t)}\right)$.
The total degree of the defining polynomial of $\mathcal{C}$ is then computed by a resultant.

The paper is organized as follows. In section 2, we recall some results of \cite{sendra1} on the implicitization of curves, and we use these results to prove Theorem \ref{main}.
Next, we present two applications of this result for constant width curves
and rotors (which is a a generalization of planar constant width curves) with polynomial support function. 
The last section is devoted to the computation of partial degrees of the defining polynomial
of algebraic surfaces which are given by a generalization of (\ref{def-curve-C}) to the euclidean space via spherical harmonics, see \cite{Groemer}.
The surfaces we consider are defined
by a bivariate polynomial which coincides with the support function when the domain enclosed by the surface is convex.
Following a procedure presented in \cite{sendra2}, we obtain the partial degrees of the implicit equation defining the
surface. The method we use allows to obtain these degrees via Maple in a reasonable time for harmonics of low order.
For these example, we have also tried a direct computation of the Gröbner basis (of the ideal generated by all polynomial vanishing on this surface) 
via Maple, but this computation could not be completed in a reasonable time.
Theses computation allow to perform a similar conjecture as Rabinowitz in the euclidean space.
The paper concludes with two sections containing the proof of the results of section 2, and a numerical code in Maple which 
was used to perform the computation of the partial degrees of the defining polynomial for surfaces.
\section{Degree of the defining polynomial of a curve of $\mathcal{E}$}
\subsection{Implicitization of a curve}
In this subsection, we recall an implicitization result of \cite{sendra1}. First, let us set some basic definitions.
We say that a property holds for {\it{almost all}} values of a parameter $t\in \mathbb{R}$ if this property is verified
for all $t \in \mathbb{R} \backslash I$, where $I$ is a finite set. If $P \in \mathbb{R}[X]$, we write $\deg P$ the {\it{degree}} of $P$, and if
$P \in \mathbb{R}[X_1,...,X_s]$, we write $\deg_{x_i} P$ the {\it{partial degree}} of $P$ with respect to $X_i$.
The {\it{total degree}} of $P \in \mathbb{R}[X_1,...,X_s]$ is the maximum of $i_1+\cdots+i_s$ of any term $X_1^{i_1} \cdots X_s^{i_s}$ of $P$.

Let $\mathcal{C}$ be a curve in $\mathcal{E}$. As $p$ is a trigonometric polynomial, $\mathcal{C}$ is algebraic,
and irreducible. Hence, it admits an irreducible defining polynomial $f\in \mathbb{R}[x,y]$ which is unique up to a multiplication
by a constant $c\not=0$, see \cite{fulton}.
An implicit equation of $\mathcal{C}$ can be obtained as follows. The composition
\begin{equation}{\label{compo}}
\mathcal P:
  \begin{array}[t]{ccc}
    \mathbb R &\to & \mathcal C\\ t &\mapsto& \mathcal \gamma(\varphi^{-1}(t)),
  \end{array}
\end{equation}
where $\varphi$ is the bijection:
\begin{equation}\
\varphi:
  \begin{array}[t]{ccc}
    ]-\pi,\pi[&\to & \mathbb R\\ \theta &\mapsto& \tan \frac \theta2,
  \end{array}
\end{equation}
provides a rational parametrization of the curve $\mathcal{C}$.
The following definition can be found in \cite{sendra1}. 
\begin{defi}
A parametrization $\mathcal{P}(t)$ is proper if and only if for almost all values of the parameter $t$,
the mapping $\mathcal{P}$ is rationally bijective.
\end{defi}
When $\mathrm{Card}(\mathcal P^{-1}(\mathcal P(t)))=1$ for almost every $t \in \mathbb{R}$, the curve $\mathcal{C}$ is proper.
The next result is a simple rephrasing of Theorem~1, Theorem~3, Theorem~5 and  Theorem~7
in~\cite{sendra1}, and gives a constructive approach to obtain the defining
  polynomial of a proper curve $\mathcal{C}$ parametrized by $\mathcal{P}(t)$.
Next, we use the notations
  $\deg\left(\frac{P_i(t)}{Q_i(t)}\right)=\max(\deg(P_i(t)),\deg(Q_i(t)))$, and $\gcd(P_i,Q_i)$ denotes the greatest common divisor (gcd) of $P_i$ and $Q_i$.
The resultant of two polynomials $H_1$ and $H_2$ is by definition the determinant of the Sylvester matrix associated to $H_1$ and $H_2$.
\begin{prop}\label{prop:radical}
  Let $\mathcal C$ be the algebraic curve over $\mathbb R$
  parametrized by $\mathcal P(t)$. Assume that for almost all
  $t\in\mathbb R$,
  \[\mbox{Card}(\mathcal P^{-1}(\mathcal P(t)))=1.\]
  Then, the defining polynomial $f(x,y)$ of $\mathcal C$ is irreducible and can
  be obtained as the resultant
  $$ f(x,y)=Res_t(H_1(t,x),H_2(t,y)),$$ where $\mathcal
  P(t)=\left(\frac{P_1(t)}{Q_1(t)},\frac{P_2(t)}{Q_2(t)} \right)$ with
  $\gcd(P_i,Q_i)=1$ and $H_1(t,x)=xP_1(t)-Q_1(t)$,
  $H_2(t,y)=yP_2(t)-Q_2(t)$.  Moreover, we have:
  $$ \deg\left(\frac{P_1(t)}{Q_1(t)}\right)=\deg_y(f) \ \ \mathrm{ and } \
  \deg\left(\frac{P_2(t)}{Q_2(t)}\right)=\deg_x(f). $$
\end{prop}
In other words, the resultant with respect to the parameter of the polynomials obtained from the parametrization
$\mathcal{P}$ by eliminating the denominator is the defining polynomial of $\mathcal{C}$.
When the number of times the curve traces $\mathcal{C}$ is greater than $1$, we have the following result, see \cite{sendra1}.
\begin{prop}{\label{non-proper}}
Let $\mathcal C$ be the algebraic curve over $\mathbb R$
  parametrized by $\mathcal P(t)$. Assume that there exists $r \in \mathbb{N}^*$ such that for almost all
  $t \in\mathbb R$,
  \[\mbox{Card}(\mathcal P^{-1}(\mathcal P(t)))=r\geq 1.\]
  Then, the defining polynomial $f(x,y)$ of $\mathcal C$ can
  be obtained as
$$
f(x,y)^r=Res_t(H_1(t,x),H_2(t,y)).
$$
Moreover, if $n=\max(\mathrm{deg}_x(f),\mathrm{deg}_y(f))$, one has for almost all $t\in \R$:
$$ r=\frac{\mathrm{deg}(\mathcal{P}(t))}{n} \ \ \mathrm{ and } \ \deg(\mathcal{P}(t))=\max\left(\deg\left(\frac{P_1(t)}{Q_1(t)}\right),\deg\left(\frac{P_2(t)}{Q_2(t)}\right)\right).$$
\end{prop}
\subsection{Main result}
The purpose of this subsection is to state our main result (Theorem \ref{main}). The proof of this result 
is divided into Lemmas \ref{index} and \ref{param}.
Lemma \ref{index} establishes the number of times $\mathcal{P}(t)$ traces $\mathcal{C}$, and lemma \ref{param}
provides a rational parametrization of $\mathcal{C}$ as in Propositions \ref{prop:radical} and \ref{non-proper}.
Let $p$ a trigonometric polynomial given by:
\begin{equation}{\label{def-p}}
p(\theta)=a_0+ \sum_{1 \leq k \leq N}(a_k \cos k \theta + b_k \sin k \theta),
\end{equation}
where $N \geq 1$, and $(a_k)_{0 \leq k \leq N}$, $(b_k)_{1 \leq k \leq N}$ are real coefficients, and let $\mathcal{C}$
the corresponding curve defined by (\ref{def-curve-C}). When $N=1$, the curve $\mathcal{C}$ defined by $p$ is a circle of radius $\frac{a_0}{2}$ and of center $(a_1,b_1)$, which has a defining polynomial of degree $2$. Without any loss of generality, we may assume $N\geq 2$ in the following.
\begin{lem}{\label{index}}Let $q$ the number of times $\gamma$ traces $\mathcal{C}$.\\
(i) If there exists $0 \leq k \leq N$ such that $a_{2k} \not=0$ or $b_{2k} \not=0$, then
$q=1$.\\
(ii) If $a_{2k}=0$ and $b_{2k}=0$ for all $k$ such that $0 \leq 2k \leq N$, then $q=2$.
\end{lem}
By using the composition (\ref{compo}), we get the following rational parametrization of the curve $\mathcal{C}$.
\begin{lem}{\label{param}} There exist two polynomials $P_1,P_2\in \mathbb{R}[t]$ satisfying
$\deg(P_1)\leq 2N+2$, $\deg(P_2)\leq 2N+1$, $\gcd(P_1,Q)=\gcd(P_2,Q)=1$, where $Q(t):=(1+t^2)^{N+1}$ and such that:
$$
\gamma(\vphi^{-1}(t))=\left(\frac{P_1(t)}{Q(t)},\frac{P_2(t)}{Q(t)} \right).
$$
\end{lem}
Combining the two previous lemmas yields to our main result.
\begin{theo}{\label{main}} Consider a trigonometric polynomial $p$ given by \eqref{def-p} with $N \geq 2$ and $a_N \not=0$ or $b_N \not=0$.
\\
(i) If there exists $0 \leq k \leq N$ such that $a_{2k} \not=0$ or $b_{2k} \not=0$, then the curve $\mathcal{C}$
has a defining polynomial of total degree $2N+2$.\\
(ii) If $a_{2k}=0=b_{2k}=0$ for all $k$ such that $0 \leq 2k \leq N$, then, the curve $\mathcal{C}$ has a defining polynomial of total degree $N+1$.
\end{theo}
\subsection{Defining polynomial of algebraic planar constant width bodies and rotors}
\begin{figure}[h!]
\begin{center}
\includegraphics[scale=0.4]{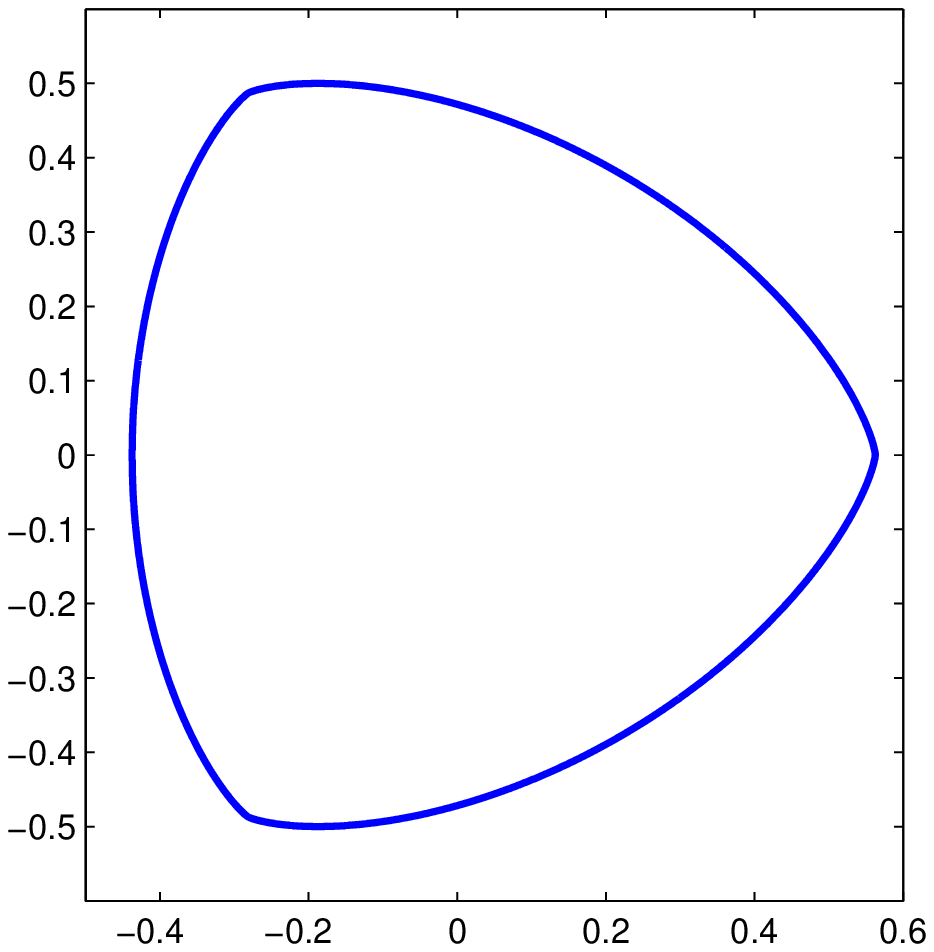}
\includegraphics[scale=0.4]{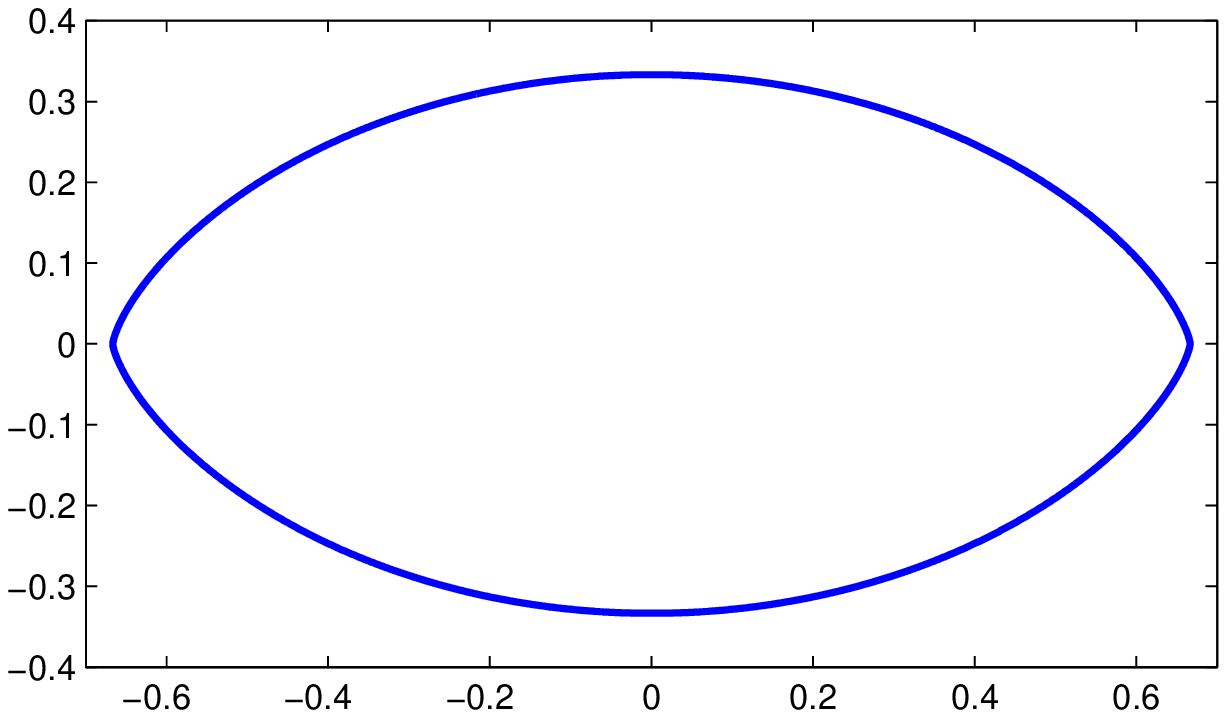}
\end{center}
\caption{{\it{Figure left}}: plot of the constant width curve defined by (\ref{c1}). {\it{Figure right}}: plot of the rotor of $P_3$ defined by (\ref{c2}). \label{fig1}}
\end{figure}

In this subsection, we present an application of Theorem \ref{main} for planar constant width curves and rotors with polynomial
support function. For future reference, let us recall that if $p$ is a given periodic function of period $2\pi$ then, the domain inside a curve $\mathcal{C}$ given by (\ref{def-curve-C}) is convex
if and only if
\begin{equation}{\label{convCourb}}
p+p'' \geq 0,
\end{equation}
in the distribution sense, see e.g. \cite{Groemer,vitale}.
\begin{defi}
We say that $\mathcal{C}$ is an algebraic curve of constant width $\alpha>0$ if and only if
$\mathcal{C}$ is described by (\ref{def-curve-C}), where $p$ is a trigonometric polynomial given by:
\begin{equation}{\label{cst-width}}
p(\theta)=\frac{\alpha}{2}+\sum_{0\leq k \leq N}(a_{2k+1} \cos (2k+1)\theta + b_{2k+1} \sin(2k+1)\theta),
\end{equation}
and that satisfies (\ref{convCourb}) .
\end{defi}
Notice that (\ref{cst-width}) implies that $p$ satisfies:
\begin{equation}{\label{largeur}}
\forall \ \theta \in [-\pi,\pi], \ p(\theta)+p(\theta+\pi)=\alpha,
\end{equation}
which is the geometrical definition of constant width curves \cite{schneider}.
The constraint $p+p'' \geq 0$ on the function $p$ is equivalent to the non-negativity of the radius of curvature
of $\mathcal{C}$ which ensures the convexity of the domain inside $\mathcal{C}$, see \cite{Bayen2},\cite{Har}.
When $p$ is a trigonometric polynomial, \eqref{convCourb} has to be understood in the classical sense.
Geometrical properties of these curves can be found in the literature, see e.g. \cite{BayenJBHU,Blaschke1,Bon,How,Yag}.
Applying Theorem \ref{main} with a polynomial function $p$ satisfying (\ref{convCourb})-(\ref{cst-width}) yields to the following result.
\begin{theo}{\label{constant-width}}
Let $\mathcal{C}$ be a non-circular algebraic curve of constant width $\alpha>0$ such that $a_{2N+1} \not=0$ or $b_{2N+1} \not=0$. Then,
the total degree of the defining polynomial $f$ of $\mathcal{C}$ is $4N+4$. Moreover, $\deg_x(f)=\deg_y(f)=4N+4$.
\end{theo}
The conjecture of Rabinowitz (see \cite{Rabi}) is then a consequence of the previous theorem.
\begin{coro}{\label{rabi-proof}}
The minimal degree of an implicit defining equation for a non-circular algebraic constant width curve is 8.
\end{coro}
The degree $8$ can be obtained by any support function $p$ satisfying given by:
\begin{equation}{\label{CWminimal}}
p(\theta)=\frac{\alpha}{2}+a_1 \cos \theta + b_1 \sin \theta + a_3 \cos 3\theta + b_3 \sin 3 \theta,
\end{equation}
where $a_3\not=0$ or $b_3 \not=0$ are small enough to ensure the convexity constraint \eqref{convCourb}.
Following Proposition \ref{prop:radical} one can
compute the defining polynomial of a curve given by \eqref{CWminimal} with $a_1=b_1=b_3=0$, $a_3=\frac{1}{16}$ and $\alpha=1$ using Maple. One has
$p(\theta)=\frac{1}{2}+\frac{1}{16} \cos 3\theta$ and $p(\theta)+p''(\theta)=\frac{1}{2}(1-\cos 3\theta) \geq 0$  (which ensures convexity of the domain inside the curve), and the parametrization of the corresponding
curve $\mathcal{C}$ becomes:
\begin{equation}{\label{c1}}
\begin{cases}
x(\theta)=\frac{1}{2}\cos \theta+\frac{1}{16}(\cos 3\theta \cos \theta+3\sin 3\theta \sin \theta),\\
y(\theta)=\frac{1}{2}\sin \theta+\frac{1}{16}(\cos 3\theta \sin \theta - 3 \sin 3\theta \cos \theta),
\end{cases}
\end{equation}
see Figure \ref{fig1}. We obtain the following defining polynomial for $\mathcal{C}$:
{\small{
\begin{align*}
&f(x,y):=-182284263-469762048y^6+1269789696y^2+33554432x^7-490733568x^4y^2+1610612736x^4y^4\\
& -134217728x^5y^2 +268435456x^6y^2+16777216x^8+6794772480y^4x+9437184x^6+4294967296y^8+317447424x^2\\
& -141557760x^5-1066991616y^2x^2+4294967296y^6x^2+2063597568y^4x^2-2133983232y^4-133373952x^4\\
& -931627008y^2x+77635584x^3-6442450944y^6x+1132462080x^3y^2-2684354560x^3y^4.
\end{align*}
}}
Similarly as constant width curves, one can define planar algebraic {\it{rotors}} as follows (see \cite{Bayen2,did-ter,Gol3,Gru}).
In the following, we say that an $n$-gon is a regular polygon with $n$ sides.
\begin{defi}
Let $n\geq 3$ and $P_n$ an n-gon with apothem $\rho>0$ (that is the radius of the inscribed circle).
We say that $\mathcal{C}$ is an algebraic rotor of $P_n$
if and only if $\mathcal{C}$ is described by (\ref{def-curve-C}), where $p$ is a trigonometric polynomial given by:
\begin{equation}{\label{rotor}}
p(\theta)=\rho+\sum_{0\leq k \leq m}(a_{kn+1} \cos (kn+1)\theta + b_{kn+1} \sin(kn+1)\theta+a_{kn-1} \cos (kn-1)\theta + b_{kn-1} \sin(kn-1)\theta),
\end{equation}
and that satisfies (\ref{convCourb}).
\end{defi}
In other words, the support function of a rotor has a Fourier expansion over harmonics of order $kn \pm 1$.
A consequence of (\ref{rotor}) is the following (see \cite{Bayen2}):
$$
\forall \ \theta \in [-\pi,\pi], \ p(\theta)-2 \cos \left(\frac{2\pi}{n}\right)p \left(\theta+\frac{2\pi}{n}\right)+ p \left(\theta+\frac{4\pi}{n}\right)=4\rho\sin^2 \left(\frac{\pi}{n}\right).
$$
For $n=4$, the previous expression simplifies into (\ref{largeur}), so that a constant width curve is a
rotor of a square. Next, we get the following result by applying Theorem \ref{main} for a function $p$ satisfying (\ref{convCourb})-(\ref{rotor}).
\begin{theo}{\label{thm-rotor}}
Let $\mathcal{C}$ be a non-circular rotor of an $n$-gon such that $a_{nm+1} \not=0$ or $b_{nm+1} \not=0$. Then,
the total degree of the defining polynomial $f$ of $\mathcal{C}$ is $2nm+4$. Moreover, $\deg_x(f)=\deg_y(f)=2nm+4$.
\end{theo}
This result coincides with Theorem \ref{constant-width} when $n=4$ (indeed, one has in this case $2N+1=4m+1$ so that $4N+4=8m+4$). Similarly
as for algebraic constant width curves, one can characterize the non-circular algebraic rotors whose defining polynomial is of minimal degree.
An immediate consequence of Theorem \ref{thm-rotor} is the following.
\begin{coro}
The minimal degree of an implicit defining equation for a non-circular algebraic rotor is $2n$.
\end{coro}
By taking $n=4$, one recovers Corollary \ref{rabi-proof}. More generally, let us now consider the set $\mathcal{R}$ of all planar algebraic rotors that are defined by (\ref{convCourb})-(\ref{rotor}) for a certain value of $n\geq 3$. We have the following result.
\begin{theo}
The minimal degree of an implicit defining equation for a non-circular algebraic curve of $\mathcal{R}$ is 6.
\end{theo}
The degree $6$ is obtained by any support function 
$$
p(\theta)=\frac{1}{2}+a_1 \cos (\theta)+ b_1 \sin(\theta)+ a_2\cos(2\theta) + b_2\sin(2\theta),
$$
where $a_2\not=0$ or $b_2\not=0$ small enough to ensure convexity constraint \eqref{convCourb}.
A rotor $\mathcal{C}$ admitting a defining polynomial of degree $6$ can be constructed as follows. Let $p(\theta)=\frac{1}{2}+\frac{1}{6}\cos(2\theta)$ (so that $p+p'' =\frac{1}{2}(1-\cos 2\theta)\geq 0$ to ensure convexity). The parametrization of the corresponding curve $\mathcal{C}$ becomes:
\begin{equation}{\label{c2}}
\begin{cases}
x(\theta) = \frac{1}{2} \cos \theta + \frac{1}{4} (\cos 2\theta \cos \theta+2\sin 2\theta \sin \theta),
\\
y(\theta) = \frac{1}{2} \sin \theta + \frac{1}{4} (\cos 2\theta \sin \theta - 2\sin 2\theta \cos \theta),
\end{cases}
\end{equation}
see Figure \ref{fig1}. The defining polynomial of $\mathcal{C}$ is:
{\small{
\begin{align*}
&f(x,y)=191102976 y^6 + 573308928 y^4  x^2+318504960 y^4  + 573308928 y^2  x^4
- 509607936 y^2  x^2  + 113246208 y^2 \\
&+191102976 x^6- 254803968 x^4+113246208 x^2  - 16777216.
\end{align*}
}}
\section{Computation of the partial degrees for parametrized surfaces}
Our aim in this section is to investigate the implicitization of several algebraic surfaces which are given by an extension of (\ref{def-curve-C}) to the euclidean space. Following the question raised by \cite{Rabi} in the planar case, we are interested in particular at studying the degree and partial degrees of the defining polynomial of algebraic constant width surfaces. Recall that a surface of constant width $\alpha>0$ is the boundary of a convex set $K\subset \R^3$ which satisfies:
$$
\forall \ \nu \in \mathcal{S}^2,  h_K(\nu)+h_K(-\nu)=\alpha,
$$
where $h_K:\mathcal{S}^2 \rightarrow \R_+$ is the support function of $K$ defined by 
\be{\label{suppFunction}}
h_K(\nu):=\sup_{x \in K} x \cdot \nu, \ \nu \in \R^n \backslash \{0\}.
\ee
Following \cite{sendra2}, we first recall how the defining polynomial $f$ and the partial degrees of $f$ 
can be obtained given a rational parametrization of a surface.
We then apply this procedure using Maple on two families of surfaces in order to find the partial degrees of the defining polynomial:
\begin{itemize}
\item For surfaces of revolution which are obtained by rotation of a planar curve given by \eqref{def-curve-C} around one axis
(not necessarily axis of symmetry of the planar curve).
\item For surfaces which are parametrized by a bivariate polynomial $h$ which is decomposed into spherical harmonics  
(and which represents the support function \eqref{suppFunction} of the set $K$ enclosed by the surface in the convex case).
\end{itemize} 
\subsection{Description of the implicitization algorithm}
All the definitions and properties are taken from \cite{sendra2}. Let us consider a rational surface over $\mathbb{R}$ given by:
$$
\mathcal{P}(t_1,t_2):=\left(\frac{p_1(t_1,t_2)}{q_1(t_1,t_2)},\frac{p_2(t_1,t_2)}{q_2(t_1,t_2)},\frac{p_3(t_1,t_2)}{q_3(t_1,t_2)} \right),
$$
where $p_i,q_i \in \mathbb{R}[t_1,t_2]$. The degree of $\mathcal{P}$ is by definition the degree of the rational map $\mathcal{P}$.
\begin{defi}
We say that the parametrization $\mathcal{P}$ is $(1,2)$ settled if and only if for all $(t_1,t_2) \in \R^2$, the vectors 
$\{ \nabla \mathcal{P}_1(t_1,t_2),\nabla \mathcal{P}_2(t_1,t_2)\}$ are linearly independent and $\mathcal{P}_3(t_1,t_2)$ is not constant.
\end{defi}
We assume that none of the projective curves obtained from $p_i(t)$ and $q_i(t)$ passes through the point at infinity $(0:1:0)$. 
We say that the {\it{general assumptions}} are verified for $\mathcal{P}$ if this property holds together with the $(1,2)$ settled parametrization.
Next, we denote by $\mathrm{lcm}(P_i,Q_i)$ the lowest common multiple of polynomials $P_i$ and $Q_i$, $\mathrm{pp}_x(P)$ 
the primitive part of polynomial $P$ (that is, the gcd of all coefficients of $\mathrm{pp}_x(P)$ is $1$), 
and $\mathrm{Content}(P)$ the Content of $P$ (the gcd of all the coefficients of $P$). Note that we have $P(x)=\mathrm{pp}_x(P)\mathrm{Content}(P)$.

Let $S(t_1,x)$ and $T(t_2,x)$ the two polynomials defined by:
$$
\begin{cases}
S(t_1,x)=\mathrm{pp}_x(\mathrm{Content}_Z(\mathrm{Res}_{t_2}(G_1,G_2+ZG_3))),\\
T(t_2,x)=\mathrm{pp}_x(\mathrm{Content}_Z(\mathrm{Res}_{t_1}(G_1,G_2+ZG_3))),
\end{cases}
$$
and for $i \in \{1,2,3\}$, $j \in \{1,2,3\}$, $i\not=j$, set:
$$
\begin{cases}
G_i(t,x_i)=p_i(t)-x_iq_i(t), \ G_4(t)=\mathrm{lcm}(q_1,q_2,q_3),\\
S_{i,j}(t_1,x)=\mathrm{pp}_x(\mathrm{Res}_{t_2}(G_i(t,x_i),G_j(t,x_j))),\\
T_{i,j}(t_2,x)=\mathrm{pp}_x(\mathrm{Res}_{t_1}(G_i(t,x_i),G_j(t,x_j))).
\end{cases}
$$
The next result is a simple rephrasing of Theorem~1, Theorem~6, Theorem~10 and Theorem~11
in~\cite{sendra2}, and gives a constructive approach to obtain the defining
  polynomial of a surface parametrized by $\mathcal{P}$.
\begin{theo} Let $\mathcal{S}$ be a rationnal affine surface defined by the irreducible polynomial $F(x)$, and let $\mathcal{P}(\ov{t})$ be a rational (1,2)-settled parametrization of $\mathcal{S}$ in reduced form. Then, there exists $r \in \mathbb{N}$ such that up to constants in $\mathbb{R}^*$,
$$
F(\ov{x})^r = \mathrm{pp}_{x_3}(h(\ov{x})),
$$
where:\\
(1) $h(\ov{x})=\mathrm{content}_{Z,W}(\mathrm{Res}_{t_2}(T_{1,2}(t_2,\ov{x}),K(t_2,Z,W,\ov{x})))$,\\
(2) $K(t_2,Z,W,\ov{x})=\mathrm{Res}_{t_1}(S_{1,2}(t_1,\ov{x}),G(t,Z,W,\ov{x}))$,\\
(3) $G(\ov{t},Z,W,\ov{x})=G_3(\ov{t},x_3)+ZG_1(\ov{t},x_1)+WG_2(\ov{t},x_2)$.\\
Moreover, we have for $i,j,k\in \{1,2,3\}$, $i<j$, $i\not=k$, $j\not=k$:
$$
deg(\mathcal{P})=\mathrm{deg}_{t_1}(S(t_1,x))=\mathrm{deg}_{t_2}(T(t_2,x)), \quad \mathrm{deg}_{x_k}(F)=\frac{\mathrm{deg}_{t_1}(S_{i,j})}{\mathrm{deg}_{t_1}(S)}=\frac{\mathrm{deg}_{t_2}(T_{i,j})}{\mathrm{deg}_{t_1}(S)}.
$$
\end{theo}
Notice that the content is taken in $\mathbb{R}(\nu)[\ov{t}]$ and that if $S_{1,2}$ and $T_{1,2}$ are square free we have $F(\ov{x})^{\mathrm{deg}(\mathcal{P})} = \mathrm{pp}_{x_3}(h(\ov{x}))$. The code which is given at the end of the paper is an implementation into Maple of the method described above.

\subsection{Computation of the partial degrees for surfaces of revolution}
We now apply the procedure above to find the partial degrees of the defining polynomial of a surface of revolution
that is obtained by rotation of a curve given by (\ref{def-curve-C}) around one axis (not necessarily of symmetry for the curve). 
The construction goes as follows.
First, consider a trigonometric polynomial $p$ and let $(x(\theta),y(\theta))$ the parametrization of a curve $\mathcal{C}$ 
given by (\ref{def-curve-C}). The surface given by:
\begin{equation}{\label{eq-surf-revol1}}
(\theta,\phi) \in ]-\pi,\pi[ \times [-\pi,\pi] \ \longmapsto \Lambda(\theta,\phi):=\begin{cases}
X(\theta,\phi)=x(\theta),\\
Y(\theta,\phi)=y(\theta)\sin\phi,\\
Z(\theta,\phi)=y(\theta)\cos \phi,
\end{cases}
\end{equation}
is by definition a surface of revolution around the axis $X$. If $p$ is the support function of a constant width curve, then 
\eqref{eq-surf-revol1} defines the boundary of a surface of constant width of revolution.

Following section 2.2, we obtain a rational parametrization of the surface parametrized by $\Lambda$ as follows:
\begin{equation}{\label{eq-surf-revol2}}
(t_1,t_2)\in \mathbb{R}^2 \longmapsto \mathcal{P}(t_1,t_2):=\Lambda(\vphi^{-1}(t_1),\vphi^{-1}(t_2))=\begin{cases}
X_1(t_1,t_2)=\frac{P_1(t_2)}{Q(t_2)},\\
Y_1(t_1,t_2)=\frac{P_2(t_2)}{Q(t_2)}\frac{2t_1}{1+t_1^2},\\
Z_1(t_1,t_2)=\frac{P_2(t_2)}{Q(t_2)}\frac{1-t_1^2}{1+t_1^2},
\end{cases}
\end{equation}
where $P_1$ and $P_2$ are given by lemma \ref{param}. We now compute the partial degrees of the defining 
polynomial $f\in \R[X,Y,Z]$ of the surface given by (\ref{eq-surf-revol2}) from the procedure described above (see the last section 
for the implementation in Maple).
Partial degrees of the defining polynomial of the surface given by (\ref{eq-surf-revol2}) are depicted for several choices of the function $p$
in Table \ref{fig2}.
From the numerical computations listed in Table \ref{fig2}, we aim at conjecturing the following result which says
that there exists a simple relation between the degree of $p$ and the total degree of the defining polynomial of $\Lambda$ (similarly as in the two-dimensional case).
\begin{conj}\label{conj:conj3d-1} Let $N \geq 2$, and $p$ a trigonometric polynomial given by (\ref{def-p}) and such that $a_N \not = 0$ or $b_N \not=0$.
Then, the partial degrees of the defining polynomial $f \in \mathbb{R}[x,y,z]$ of the surface of $\mathbb{R}^3$ defined by (\ref{eq-surf-revol1}) satisfy
$$
\mathrm{deg}_{x} f=\mathrm{deg}_{y} f=\mathrm{deg}_{z} f = 4N+4,
$$
and the total degree of $f$ is $4N+4$. In particular a constant width surface of revolution has a defining polynomial of degree $4N+4$ provided that 
its support function $p$ is a trigonometric polynomial of degree $N$.
\end{conj}
\begin{table}
\begin{center}
\begin{tabular}{|c|c|c|c|c|c|}
\hline
$p(\theta)$ & symmetric & $\mathrm{deg}(\mathcal{P})$ & $\frac{\mathrm{deg}_{x} f}{\mathrm{deg}(\mathcal{P})}$ & $\frac{\mathrm{deg}_{y}f}{\mathrm{deg}(\mathcal{P})}$ & $\frac{\mathrm{deg}_{z}f}{\mathrm{deg}(\mathcal{P})}$\\
\hline
$\cos 2\theta$ & yes & 2 & 6 & 6 & 6 \\ \hline
$\cos 3\theta$ &  yes & 4 & 4 & 4 & 4 \\ \hline
$\cos 4\theta$ & yes & 2 & 10 & 10 & 10 \\\hline
$\cos 5 \theta$ & yes & 4 & 6 & 6 & 6 \\\hline
$\frac{1}{2}+\frac{1}{32}\cos 3\theta$ & yes & 2 & 8 & 8 & 8 \\\hline
$\frac{1}{2}+\frac{1}{32}\cos 3\theta+\sin3\theta$ & no & 1 & 16 & 16 & 16 \\\hline
$\frac{1}{2}+\sin 3\theta$ & no & 1 & 16 & 16 & 16 \\ \hline
$\cos 4\theta +\sin 3\theta$ & no & 1 & 20 & 20 & 20 \\ \hline
$\cos 4\theta + \sin 5\theta$ & no & 1 & 24 & 24 & 24 \\
\hline
\end{tabular}
\end{center}
\caption{Partial degrees for surfaces of revolution parametrized by (\ref{eq-surf-revol2}). \label{fig2}}
\end{table}
\subsection{Computation of the partial degrees for surfaces described by spherical harmonics}
In this subsection, we consider surfaces of $\R^3$ parametrized by:
\begin{equation}{\label{param1}}
(\theta,\phi) \in ]0,\pi[ \times [-\pi,\pi] \longmapsto \Upsilon(\theta,\phi) = h u + h_{\theta} v + \frac{1}{\sin \theta} h_{\phi} w,
\end{equation}
where $(u,v,w)$ is the orthonormal basis:
\begin{center}
$u=$
\begin{tabular}{|c}
$\sin \theta \cos \phi$ \\
$\sin \theta \sin \phi$ \\
$\cos \theta$
\end{tabular}
, \
$v=$
\begin{tabular}{|c}
$\cos \theta \cos \phi$ \\
$\cos \theta \sin \phi$ \\
$-\sin \theta$
\end{tabular}
, \
$w=$
\begin{tabular}{|c}
$-\sin \phi$ \\
$\cos \phi$ \\
$0$
\end{tabular}
,
\end{center}
and $h$ is decomposed into a finite sum of spherical harmonics (see \cite{Groemer,zbynek}):
\begin{equation}{\label{dev1}}
h(\theta,\phi)=a_{0,0} + \sum_{2 \leq l \leq N} \sum_{0 \leq m \leq l}P_l^m(\cos \theta)(a_{l,m} \cos m\phi + b_{l,m} \sin m\phi).
\end{equation}
The mapping $u \longmapsto P_l^m(u)$ represents the Legendre associate function (see \cite{Groemer}):
$$
P_l^m(u)= \frac{1}{2^l l!}(-1)^{m+l} (1-u^2)^{\frac{m}{2}}\frac{d^{l+m}}{du^{l+m}}(1-u^2)^l, \ l \geq 0, \ 0 \leq m \leq l.
$$
It is well-known that
given a non-empty compact convex set $K$ of $\mathbb{R}^3$ with strictly convex boundary, then its boundary
is described by (\ref{param1}), where $h$ is the support function of $K$, see \cite{Bon}. 
Conversely, given a function $\sigma:\mathcal{S}^2 \rightarrow \mathbb{R}$ of class $C^1$,
one can define a function $h:\mathbb{R}^2 \rightarrow \mathbb{R}$ by $(\theta,\phi) \longmapsto h(\theta,\phi)=\sigma \circ u(\theta,\phi)$ which is periodic in $(\theta,\phi)$ of period $2\pi$. The corresponding surface $\mathcal{S}$ defined by (\ref{param1}) is not necessarily the boundary of a convex set, see \cite{Bayen4}. When $h$ is given by (\ref{dev1}), (\ref{param1}) defines the boundary of a convex set provided
that coefficients $(a_{l,m},b_{l,m})$ are small enough (see \cite{did-ter} for a study of this problem in the two-dimensional framework).

When $a_{l,2k}=b_{l,2k}=0$ for $l\geq 1$, one recovers by (\ref{param1}) algebraic surfaces of constant width (see \cite{Bayen4,BayenJBHU} for a study of constant width surfaces or {\it{spheroforms}}). 
In the particular case where $h$ does not depend on $\phi$, the surface given by \eqref{param1} is a surface of constant width of revolution.

Similarly as in the two-dimensional case, one can study the question of implicitization of algebraic surfaces of constant width, and the total degree of the defining polynomial of a non-circular algebraic constant width surface.
The exact determination of these degrees seems a difficult question, but we can obtain them for low degrees of $h$ following \cite{sendra2}.

If $N \leq 1$, the formula (\ref{dev1}) writes
$h(\theta,\phi)=a_{00}P_0^0(\cos \theta)+P_1^0(\cos \theta)+P_1^1(\cos \theta)(a_{11}\cos\phi+b_{11} \sin \phi)$,
where $P_0^0(\cos \theta)=1$, $P_1^0(\cos \theta)=1=\cos \theta$, $P_1^1(\cos \theta)=-\sin \theta$. Thus, (\ref{param1}) becomes
\begin{center}
$\Upsilon(\theta,\phi)=$
\begin{tabular}{|c}
$-a_{11}+a_{00}\sin \theta \cos \phi$, \\
$-b_{11}+a_{00}\sin \theta \sin \phi$, \\
$a_{10}+a_{00}\cos \theta$,
\end{tabular}
\end{center}
which represents the sphere of radius $a_{00}$ and of center $(-a_{11},-b_{11},a_{10})$, therefore we may assume $N \geq 2$ in the numerical computations presented in Tabular \ref{fig3}.
In the following, we consider a rational parametrization of (\ref{param1}) which is obtained as follows:
$$
(t_1,t_2)\in \mathbb{R}^2 \longmapsto \mathcal{P}(t_1,t_2):=\Lambda(\vphi^{-1}(t_2),\vphi^{-1}(t_1)),
$$
and we compute the partial degrees of the defining polynomial of $\mathcal{P}$ for different choices of $h$ given by (\ref{dev1}).
In tabular \ref{fig3}, the partial degrees of the defining polynomial of a surface given by (\ref{param1}) are depicted for several choices
of the spherical harmonics. We observe that there exist spherical harmonics such that $\mathrm{deg}_{x} f \not=\mathrm{deg}_{z} f$,
whereas for surfaces of revolution, we have obtained $\mathrm{deg}_{x} f =\mathrm{deg}_{y} f = \mathrm{deg}_{z} f$ for low degrees.
\begin{table}
\begin{center}
\begin{tabular}{|c|c|c|c|c|}
\hline
$h(\theta,\phi)$ & $\mathrm{deg}(\mathcal{P})$ & $\frac{\mathrm{deg}_{x} f}{\mathrm{deg}(\mathcal{P})}$ & $\frac{\mathrm{deg}_{y}f}{\mathrm{deg}(\mathcal{P})}$ & $\frac{\mathrm{deg}_{z}f}{\mathrm{deg}(\mathcal{P})}$\\\hline
$P_2^0(\cos \theta)$ & 2& 6& 6& 6 \\ \hline
$P_2^1(\cos \theta)(\cos \phi + \sin \phi)$ & 2 & 10 & 10 & 8 \\ \hline
$P_2^2(\cos \theta)(\cos 2\phi + \sin 2\phi)$ & 2& 10 & 10 & 10 \\ \hline
$P_3^0(\cos \theta)$ & 4 & 4 & 4 & 4 \\ \hline
$1+ P_3^0(\cos \theta)$ & 2 & 8 & 8 & 8 \\ \hline
$P_3^1(\cos \theta)(\cos \phi + \sin \phi)$ & 4& 6 & 6 & 6 \\ \hline
$1+P_3^1(\cos \theta)( \cos \phi + \sin \phi)$ & 2 & 14 & 14& 12\\ \hline
$P_3^2(\cos \theta)(\cos 2\phi + \sin 2\phi)$ & 4 & 10 & 10 & 8 \\ \hline
$P_3^3(\cos \theta)(\cos 3\phi + \sin 3\phi)$ & 4 & 10 & 10 & 10 \\ \hline
$1+P_3^3(\cos \theta)(\cos 3\phi + \sin 3\phi)$ & 2 & 20 & 20 & 20 \\ \hline
$P_4^0(\cos \theta)$ & 2 & 10 & 10 & 10 \\ \hline
$P_4^1(\cos \theta)(\cos \phi + \sin \phi)$ & 2 & 18 & 18 & 16 \\ \hline
$P_4^2(\cos \theta)(\cos 2\phi + \sin 2\phi)$ & 2 & 30 & 30 & 26 \\ \hline
$P_4^3(\cos \theta)(\cos 3\phi + \sin 3\phi)$ & 2 & 34& 34 & 28\\ \hline
$P_4^4(\cos \theta)( \cos 4\phi + \sin 4\phi)$ & 2 & 34& 34& 34\\ \hline
\end{tabular}
\end{center}
\caption{Partial degrees for surfaces parametrized by (\ref{param1})-(\ref{dev1}) for different
values of $(l,m)$. \label{fig3}}
\end{table}
\subsection{Discussion}
This paper has tackled the question of implicitization of a certain family of algebraic curves and surfaces.  
In the planar case and when the curve is given by a trigonometric polynomial $p$, 
we have obtained the degree of the defining polynomial in term of the degree of $p$. In addition, 
we could apply this result in the particular case of constant width curves and rotors. 
Notice that the convexity plays no role in this study.

In the three-dimensional case, the surface is given via a trigonometric polynomial $h$ which is decomposed into spherical harmonics. 
Applying the same procedure in this case for determining the degree of the defining polynomial seems more delicate 
in view of the implicitization algorithm. Thanks to the theory of elimination, we have determined the degree of this polynomial for low degree of $h$. 
This computation could not be obtained via Gröbner basis in a reasonable amount of time.
We can expect proving Conjecture~\ref{conj:conj3d-1} 
(by using the implicitization algorithm) as the parametrization 
of a surface of revolution only slightly differs from the planar parametrization.
However, in the more general case, finding the degree of the defining polynomial of an algebraic constant width surface seems a difficult question.
\section{Proof of the main results}
{\bf{Proof of lemma \ref{index}}}. (i)
By derivating (\ref{def-curve-C}), one has $\gamma'(\theta)=\rho(\theta)v_{\theta}$, where
$\rho(\theta)=p(\theta)+p''(\theta)$, and $v_{\theta}=(-\sin \theta,\cos \theta)$.
Let $r=Card(\mathcal{P}^{-1}(\mathcal{P}(\alpha))$, and assume that $r \geq 2$. The equation $\rho(\theta)=0$ has
a finite number of solutions on $[0,2\pi]$ as $p$ is polynomial, and let $V$ this set. As the polynomial $\rho_1(\theta):=\rho(\theta)+\rho(\theta+\pi)$
has a non-zero coefficient, the polynomial $\rho_1$ has a finite number of zeroes on $[0,2\pi]$, and let $W$ this set.
Both sets $\gamma(V)$ and $\gamma(W)$
consist of a finite number of points of $\mathcal{C}$.
Let $M$ be a point on $\mathcal{C}\backslash (\gamma(V)\cup \gamma(W))$
which is regular (recall that the number of singular points of $\mathcal{C}$ is finite as $\mathcal{C}$ is algebraic).
By definition of $r$, there exists $(t_1,t_2)\in \mathbb{R}^2$ such that $t_1\not=t_2$, and $\mathcal{P}(t_1)=\mathcal{P}(t_2)$.
If $\theta_i=\vphi^{-1}(t_i)$, $i=1,2$, we have $M=\gamma(\theta_1)=\gamma(\theta_2)$, and $\theta_1 \not=\theta_2$ as $\vphi$ is a bijection.
We may assume $\theta_1 < \theta_2$. Now we have $\gamma'(\theta_1)=\rho(\theta_1)v_{\theta_1} \not=0$, and $\gamma'(\theta_2)=\rho(\theta_2)v_{\theta_2} \not=0$. If $\gamma'(\theta_1)=\gamma'(\theta_2)$, then we must have $\theta_2=\theta_1+\pi$ and $\rho(\theta_2)=-\rho(\theta_1)$. Hence, $\rho(\theta_1+\pi)+\rho(\theta_1)=0$ which is a contradiction. Thus, we necessarily have 
$\gamma'(\theta_1)\not=\gamma'(\theta_2)$,
and we obtain a contradiction as $M$ is a regular point of $\mathcal{C}$.
\\
(ii) Using that $p$ has only even coefficients, one gets immediately $\gamma(\theta+\pi)=\gamma(\theta)$ for all $\theta \in [0,2\pi]$, hence,
$r=Card(\mathcal{P}^{-1}(\mathcal{P}(\alpha))$ is greater than $2$. Assume that the number of times $\gamma$ traces $\mathcal{C}$ on the
interval $[0,\pi]$ is strictly greater than $2$, and let $M$ be a regular point of $\mathcal{C}$.
There exists $0<\theta_1<\theta_2\leq \pi$ such that $M=\gamma(\theta_1)=\gamma(\theta_2)$. But, one has $\gamma'(\theta_1)=\gamma'(\theta_2)$ if and only if $\theta_1=\theta_2$ or $\theta_2=\theta_1+\pi$. Both conditions are not possible, hence, $\gamma'(\theta_1)\not= \gamma'(\theta_2)$, which
contradicts that $M$ is regular. Hence, $r=2$.
\hfill$\Box$
\\
\\
{\bf{Proof of lemma \ref{param}}}.
Let us compute a rational parametrization of the curve $\mathcal{C}$. Recall that Chebychev's polynomial of first and second kind $T_n$ and $U_n$ are given by:
$$\begin{cases} T_n(x)=\sum_{0 \leq k \leq \lfloor \frac{n}{2} \rfloor} C_n^{2k} (-1)^k
  x^{n-2k}(1-x^2)^k, \\ U_n(x)=\sum_{0 \leq k \leq \lfloor \frac{n}{2}\rfloor} C_{n+1}^{2k+1}
  (-1)^k x^{n-2k}(1-x^2)^k,
\end{cases}
$$
so that (recall $t=\tan\left(\frac{\theta}{2}\right)$):
\begin{align*}
  \cos(n\theta)&=\frac{C_n(t)}{(1+t^2)^n}& \mbox{with } & \ \ C_n(t)=\sum_{0
    \leq k \leq \lfloor \frac{n}{2} \rfloor} C_n^{2k} (-1)^k 2^{2k}t^{2k}(1-t^2)^{n-2k}, \\
  \sin(n\theta)&=\frac{S_{n}(t)}{(1+t^2)^n}&\mbox{with } & \ \ S_n(t)=\sum_{0
    \leq k \leq \lfloor \frac{n}{2} \rfloor} C_{n}^{2k+1} (-1)^k 2^{2k+1}t^{2k+1}(1-t^2)^{n-1-2k}.
\end{align*}
The fractions are irreducibles as $C_n(\pm i)=2^{2n-1}$ and $S_n(\pm
i)=\pm i 2^{2n-1}$. Moreover, polynomials $C_n$ and $S_n$ are of degree $2n$.
The parametrization (\ref{def-curve-C}) can be written:
$$
  x(\theta)=a_0\cos \theta+\sum_{k=1}^N \left(a_{k} \cos\theta\cos (k
  \theta)+b_{k} \cos\theta\sin (k \theta)+ ka_{k} \sin\theta\sin (k
  \theta)-kb_{k} \sin\theta\cos (k \theta)\right),
$$
$$
  y(\theta)=a_0\sin \theta+\sum_{k=1}^N \left(a_{k} \sin\theta\cos (k
  \theta)+b_{k} \sin\theta\sin (k \theta)- ka_{k} \cos\theta\sin (k
  \theta)+kb_{k} \cos\theta\cos (k \theta)\right).
$$
Let us define $(\tilde{x}(t),\tilde{y}(t)):=\gamma(\vphi^{-1}(t))$. The parametrization of $\mathcal{C}$ rewrites in $t$ as 
  $$
  \label{eq:param}
  \begin{cases}
 \tilde{x}(t)& =\frac{P_1(t)}{(1+t^2)^{N+1}},\\
 \tilde{y}(t)& =\frac{P_2(t)}{(1+t^2)^{N+1}},
  \end{cases}
$$
with
\begin{align*}
  P_1(t)&= a_0(1-t^2)(1+t^2)^N+\sum_{k=1}^N(1-t^2)(1+t^2)^{N-k}\left(a_{k} C_{k}(t)+b_{k} S_{k}(t)\right)\\
&~~+ \sum_{k=1}^N2t(1+t^2)^{N-k}(ka_{k} S_{k}(t)-kb_{k} C_{k}(t)),\\
  P_2(t)&= 2a_0t(1+t^2)^N+\sum_{k=1}^N 2t(1+t^2)^{N-k}\left(a_{k} C_{k}(t)+b_{k} S_{k}(t)\right)\\
&~~+ \sum_{k=1}^N(1-t^2)(1+t^2)^{N-k}(-ka_{k} S_{k}(t)+kb_{k} C_{k}(t)).
\end{align*}
By the expression above, we get:
$$
\begin{cases}
P_1(\pm i)=2^{2N}(a_N\pm i b_N)(1\mp N)\neq 0,\\
P_2(\pm i)=2^{2n}(b_N\mp i a_N)(N-1)\neq 0,
\end{cases}
$$
as $N>1$. Hence, the fractions are in reduced form. Moreover, the degree of $P_1(t)$ satisfies $\deg P_1 \leq 2N+2$ and the degree of $P_2(t)$ is such that $\deg P_2 \leq 2N+2$ which proves the lemma.\hfill$\Box$
\\
\\
{\bf{Proof of theorem \ref{main}}}. (i) In this case, the number of times $\gamma$ traces $\mathcal{C}$ is $1$ by Lemma \ref{index}, and we may apply Proposition \ref{prop:radical}. The resultant can be obtained by calculating the determinant of the
  Bezout matrix for the polynomials $H_1(t,x)$ and $H_2(t,x)$, defined
  by:
  \[B=(b_{i,j})_{0 \leq i,j\leq 2N+2} \ \mbox{ where } \ \frac{H_1(t_1,x)H_2(t_2,y)-H_1(t_2,x)H_2(t_1,y)}{t_1-t_2}=\sum_{0 \leq i,j \leq 2N+2}b_{i,j}t_1^it_2^j.\]
Thus, the total degree of $f$ is less than $2N+2$.
  Moreover, we have $\deg\left(\frac{P_j(t)}{(1+t^2)^{N+1}}\right)=2N+2$ for
  $j=1,2$.

(ii) In this case, the number of times $\gamma$ traces $\mathcal{C}$ is $2$ by Lemma \ref{index}, and we may apply Proposition \ref{non-proper}.
Similarly as for (i), we obtain that the degree of $f$ is less or equal than $2N+2$. But, one has $r=2=\frac{\deg(\mathcal{P}(t))}{n}=\frac{2N+2}{n}$,
where $n=\max(\mathrm{deg}_x(f),\mathrm{deg}_y(f))$, hence $n=N+1$.\hfill$\Box$
\section{Maple Code}
Hereafter, we give the main procedure {\texttt{verify}} that 
we have implemented into Maple to compute the partial degrees of the defining polynomial of a surface given in parametrization form. 
It makes use of another procedure {\texttt{isinGeneralAss}} which ensures that the general assumptions are verified.

\noindent
\small{
\begin{verbatim}
verif:=proc(Pp,t1t2:=[rand(10)(),rand(10)()])
local P,Z,W,tosubs,tosubs2,G,st,Rest2,res2,res3,degP,degP1,degS12,degS13,degS23,degS12b,degS13b,degS23b,
    gcd12,gcd13,gcd23,t1t2bis;
    P:=Pp;
#verification of the point at infinity
    if not isinGeneralAss(P,t[1],t[2],[0,1,0]) then
        P:=subs(t[1]=t[1]+t[2],P);
        print("change of variable t1->t1+t2");
    fi;
    if isinGeneralAss(P,t[1],t[2],[0,1,0]) then print("point [0,1,0] ok");
    else print("Problem, point [0,1,0] is not ok");fi;
# verification (1,2)-settled
    if is12settled(P[1],P[2],t[1],t[2]) then print("(1,2)-settled ok");
    else print("Problem, (1,2)-settled is not ok");fi;
# Computation of the Gi
    G:=[seq(numer(-x[i]+P[i]),i=1..3),mul(u[1],u=factors(mul(denom(P[i]),i=1..3))[2])];
# evaluation points:
    t1t2bis:=[rand(11..100)(),rand(1..10)()]:
    print("evaluation point"=t1t2,t1t2bis);
    tosubs:=[seq(x[i]=subs(t[1]=t1t2[1],t[2]=t1t2[2],P[i]),i=1..3)];
    tosubs2:=[seq(x[i]=subs(t[1]=t1t2bis[1],t[2]=t1t2bis[2],P[i]),i=1..3)];
    Rest2:=resultant(subs(tosubs,G[1]),subs(tosubs,G[2]+Z*G[3]),t[2]):
    res2:=resultant(subs(tosubs,G[1]),subs(tosubs,G[2]+Z*G[3]+W*G[4]),t[2]):
    res3:=resultant(subs(tosubs2,G[1]),subs(tosubs2,G[2]+Z*G[3]),t[2]):
    degP1:=(numer(factor(normal(content(Rest2,Z)/content(res3,Z)))));
    print(primpart(degP1,t[1]));
    degP1:=degree(degP1);
    degP:=degree(normal(content(Rest2,Z)/content(res2,[Z,W])),t[1]);
    if degP1<>degP then print("not the same degree",degP1,degP);fi;
    degP:=min(degP1,degP);
    degS12:=resultant(subs(tosubs,G[1]),subs(tosubs,G[2]),t[2]);
    degS13:=resultant(subs(tosubs,G[1]),subs(tosubs,G[3]),t[2]);
    degS23:=resultant(subs(tosubs,G[2]),subs(tosubs,G[3]),t[2]);
    degS12b:=resultant(subs(tosubs2,G[1]),subs(tosubs2,G[2]),t[2]);
    degS13b:=resultant(subs(tosubs2,G[1]),subs(tosubs2,G[3]),t[2]);
    degS23b:=resultant(subs(tosubs2,G[2]),subs(tosubs2,G[3]),t[2]);
    print(map(degree,[degS12,degS12b,degS13,degS13b,degS23,degS23b]));
    gcd12:=gcd(degS12,degS12b);print("gcd12",gcd12);
    gcd13:=gcd(degS13,degS13b);print("gcd13",gcd13);
    gcd23:=gcd(degS23,degS23b);print("gcd23",gcd23);
    if degree(gcd23)>0 then print("pgcd 23"=gcd23);degS23:=normal(degS23/gcd23);fi;
    if degree(gcd12)>0 then print("pgcd 12"=gcd12);degS12:=normal(degS12/gcd12);fi;
    if degree(gcd13)>0 then print("pgcd 13"=gcd13);degS13:=normal(degS13/gcd13);fi;
    degree(normal(degS13/gcd13))/degP,degree(normal(degS12/gcd12))/degP,t1t2];
return [degP,degree(degS23)/degP,degree(degS13)/degP,degree(degS12)/degP,t1t2];
end;
\end{verbatim}
}
\section{Acknowledgements}
The authors would like to thank Jean-Baptiste Hiriart-Urruty for fruitful discussions and for indicating us the problem. 

\end{document}